\documentclass[11pt]{amsart}

\usepackage{amsmath}
\usepackage{amssymb}
\usepackage{graphicx}
\usepackage{hyperref}
\usepackage{color}

\newtheorem{theorem}{Theorem}[section]
\newtheorem{corollary}[theorem]{Corollary}
\newtheorem{lemma}[theorem]{Lemma}

\theoremstyle{definition}
\newtheorem{definition}[theorem]{Definition}

\newtheorem{remark}[theorem]{Remark}

\newcommand{\la}{\lambda }

\newcommand{\lan}{\langle }
\newcommand{\ran}{\rangle }

\newcommand{\wi}{\widetilde}
\newcommand{\rh}{\rightharpoonup}

\newcommand{\vp}{\varphi}
\newcommand {\A}{A:X\supset D(A)\to 2^{X^*}}

\numberwithin{equation}{section}

\title[General Duality Mappings and Yosida Approximants ]{Continuity of the Yosida Approximants  Corresponding to General Duality Mappings}

\author[D. R. Adhikari]{Dhruba R. Adhikari}

\address[D. R. Adhikari]{Kennesaw State University,  1100 South Marietta Pkwy, Marietta, GA 30060, USA}
\email{{\tt dadhikar@kennesaw.edu}}

\keywords{Maximal monotone operators,	duality mappings,  Yosida approximants and resolvents,  Browder degree}

\subjclass[2010]{Primary 46B20,   Secon
dary 47H05}

\begin{document}

\begin{abstract}
Let $X$ be a real locally uniformly convex Banach space and  $X^*$ be the dual space of $X$.
Let $\vp:\mathbf R_+\to \mathbf R_+$ be a strictly increasing and continuous function  such that $\vp(0) = 0$, $\vp(r) \to \infty$ as $r\to\infty$,    and let $J_\vp$ be the duality mapping corresponding to  $\vp$. 		 We will prove that for every $R>0$ and every $x_0\in X$ there exists a nondecreasing function $\psi  = \psi (R, x_0) :\mathbf R_+\to \mathbf R_+$  such that $\psi(0) = 0$, $\psi(r)>0$ for $r>0$, and 
$\lan x^*- x_0^*, x-x_0\ran \ge \psi(\|x-x_0\|) \|x-x_0\|$
for all $x$ satisfying $\|x-x_0\|\le R$ and all $x^*\in J_\vp x$ and $x_0^*\in J_\vp x_0.$ This result extends the previous results of Pr\"{u}{\ss}   and Kartsatos who studied the normalized duality mapping $J$ (with $\vp(r)=r$)  for  uniformly convex  and locally uniformly Banach spaces,  respectively.   As an application of the above result,  we give a concise proof of the continuity of the Yosida approximants $A_\la^\vp$ and resolvents  $J_\la^\vp$ of a maximal monotone operator  $\A$  on  $(0, \infty) \times X$ for an arbitrary $\vp$ when $X$ is  reflexive and both $X$ and $X^*$ are locally uniformly convex.  In addition, we  discuss  pseudomonotone homotopy  of the Yosida approximants $A_\la^\vp$ with reference to the Browder degree.
\end{abstract}

\maketitle


\section{Introduction and Preliminaries}\label{S1}

Let $X$ be a real  Banach space with $X^*$ its dual space. The symbol $2^{X^*}$ denotes the collection of all subsets of $X^*$.  In what follows, $\|\cdot\|$ denotes
the norms of both  $X$ and $X^*$, and an underlying space will be understood from the context in which the norm symbol is used.  The pairing $\langle x^*,x\rangle$ denotes the value of the
functional $x^*\in X^*$ at $x\in X$.
The symbols $\partial D,$  $\overline D$ and $ \operatorname{co}D$ denote
the  boundary, closure, and convex hull of the set $D,$  respectively.
The symbol $B_R(x_0)$ denotes the open ball of radius $R>0$ with center at
$x_0$ in a Banach space. 	The symbols $\mathbf{R}$ and $\mathbf{R}_+$ denote $ (-\infty, \infty)$ and
$[0,\infty)$, respectively.  

We recall that a real Banach space $X$ is \textit{locally uniformly convex} if for every $x_0\in X$ with $\|x_0\|= 1$ and for every $\varepsilon \in (0, 2]$,  there exists a number $\delta = \delta (\varepsilon, x_0)$ such that  for all $x\in X$ with $\|x\| =1$ and
$\|x- x_0\| \ge \varepsilon$, we have $ \|x+x_0\| \le 2 (1- \delta).$ A Banach space $X$ is \textit{strictly convex }if  for all $x, y \in X$ with $ x\ne y$  and
$\|x\| = \|y\| = 1,$ we have $ \|x+y \| <2.$ A Banach space $X$ is \textit{smooth} if for every $x\ne 0$ in $X$ there exists a unique $x^*\in X^*$ such that $\|x^*\| = 1$ and $\lan x^*, x\ran = \|x\|.$

If $\{x_n\}$ is a sequence in $X$, we denote its strong convergence to
$x_0$ in $X$ by $ x_n\to x_0$   and its weak convergence to $x_0$ in $X$
by $x_n\rightharpoonup x_0.$  An operator $A : X\supset D(A)\to X^*$ is said
to be \textit{bounded} if it maps bounded subsets of the domain $D(A)$ onto
bounded subsets of  $X^*$;
\textit{demicontinuous} at $x\in D(A)$ if $x_n\in D(A)$ and $x_0\in D(A)$, $x_n\to x_0\in X$   imply $Ax_n\rh Ax_0$  in $X^*$;  
\textit{hemicontinuous} at $x\in D(A)$ if $t_n\to 0^+$  and $x +t_ny\in D(A)$ for $y\in X$   imply $A(x +t_ny)\rh  Ax $ in $X^*;$ and 
of  \textit{type $(S_+)$ }if for every sequence $\{x_n\}\subset D(A)$, $x_0\in X$ with $x_n\rh  x_0$ in $X$  such that $ \limsup_{n\to\infty} \lan Ax_n, x_n - x_0\ran \le 0,$ we have  $x_n\to x_0$ in $X$.

%

For a multivalued operator $A$ from $X$ to $X^*$, we write
$A:X\supset D(A) \to 2^{X^*}$, where $D(A)=\{x\in X: Ax\neq\emptyset\}$
is the effective domain of $A$. 
We denote by ${\rm Gr}(A)$ the graph of $A$, i.e.,
${\rm Gr}(A)=\{(x,y): x\in D(A), y\in Ax\}$. In what follows,  $Ax$  and $A(x)$ will be used interchangeably.

A \textit{gauge} function  is a strictly  increasing continuous function $\varphi: \mathbf{R}_+ \to \mathbf{R}_+$ with $\varphi(0) = 0$ and $\varphi(r) \to \infty$ as $r\to\infty$.  The  duality
mapping of $X$ corresponding to a gauge function $\varphi $ is the mapping  $J_\varphi:X \supset D(J_\varphi )\to 2^{X^*}$  defined by
$$
J_\vp x = \{x^*\in X^* : \langle x^*, x\rangle
= \varphi(\|x\|)\|x\|, \; \|x^*\| = \varphi(\|x\|)\}, \quad x\in X.
$$
The Hahn-Banach theorem ensures that $D(J_\varphi) = X,$ and therefore  $J_\varphi :X\to 2^{X^*}$
is a multivalued mapping in general. The  duality mapping corresponding to the gauge function $\varphi (r) = r$ is called the \textit{normalized duality} mapping,  denoted simply by $J$. 
The  duality mapping $J_\varphi $  satisfies
\begin{equation}\label{1.1}
	J_\varphi x  = \frac{\vp(\|x\|)}{\|x\|} Jx, \quad  x\in X.
\end{equation}
When $X$ and $X^*$ are reflexive and strictly convex,  the inverse $J^{-1}_\vp$ of the single-valued $J_\vp$ is the duality mapping $J_{\vp^{-1}}$ of $X^*.$ 
Further properties of duality mappings and geometry of Banach spaces  can be found, for example,  in Alber and Ryazantseva~\cite{Alber} and Cioranescu~\cite{Cioranescu}.

An operator $A:X\supset D(A)\to 2^{X^*}$
is said to be \textit{monotone} if for all $x, y\in D(A)$ and all
$u\in Ax,v\in Ay$ we have
$$
\langle u - v,x-y\rangle \ge 0.
$$
A monotone operator $A$ is said to be \textit{maximal monotone} if ${\rm Gr}(A)$ is
maximal in $X\times X^*$, when $X\times X^*$ is partially ordered by 
set inclusion. Suppose $X$ and $X^*$ are reflexive and strictly convex. A monotone operator
$\A$ is maximal if and only if $R(A+\lambda J) = X^*$ for all $\lambda\in (0,\infty)$.
With this setting of $X$ and $X^*$ and a maximal monotone $A$, 	its  Yosida approximant $A_t\equiv(A^{-1}+tJ^{-1})^{-1}:X\to X^*$, $t>0$,    is bounded, demicontinuous,
maximal monotone and such that $A_tx\rightharpoonup A^{\{0\}}x$ as $t\to 0^+$ for every
$x\in D(A)$, where $A^{\{0\}}x$ denotes the element $y^*\in Ax$ of minimum norm,
i.e., $\|A^{\{0\}}x\|=\inf\{\|y^*\|: y^*\in Ax\}$.
In our setting, this infimum is always attained and
$D(A^{\{0\}})=D(A)$. Also, $A_tx\in AJ_tx$, where $J_t \equiv I-tJ^{-1}A_t:X\to X$
and satisfies $\lim_{t\to 0}J_tx=x$ for all $x\in\overline{\operatorname{co}D(A)}.$ 
The operators $A_t$ and $J_t$ were introduced by Br\'ezis, Crandall and
Pazy in \cite{BCP}. For their basic properties, we refer the reader
to \cite{BCP} as well as Pascali and Sburlan \cite[pp. 128-130]{PS}. 	For further information on functional analytic tools used herein,
the reader is referred to Barbu \cite{BA}, Browder \cite{Browder1976}, Simons \cite{Simons},
Skrypnik \cite{Skrypnik1986,Skrypnik1994}, and Zeidler \cite{Zeidler1}.

The following theorem gives a necessary and sufficient condition for a monotone operator to be maximal monotone  in terms of an arbitrary  gauge function.   Although  the theorem is expected to hold for any $\vp$, to  our knowledge,  its proof can  be found only for   $\vp(r) = r^{p-1}$,  $p\in (1, \infty),$   in the literature (e.g. \cite{Alber,BA, PS}). Since the theorem will be applied when we define  general Yosida approximants and resolvents in  Section~\ref{S3},  we give  its complete proof here.

\begin{theorem}\label{Th1}
	Let $X$ be a   reflexive Banach space and both $X$ and $X^*$ strictly convex. 	A monotone operator $\A$ is maximal if and only if $R(A+\lambda J_\vp) = X^*$ for all $\lambda\in (0,\infty)$ and all  gauge functions $\vp.$ 
\end{theorem}

\begin{proof}
	Suppose that  $A$ is monotone.  The duality mapping $J_\vp :X\to X^*$ is single-valued becuase $X^*$ is strictly convex ($X$ is smooth).  Fix a $\la>0$ and a gauge function $\vp$, and suppose that $R(A+\lambda J_\vp) = X^*.$  Let $(x_0, y_0) \in X\times X^*$ satisfy 
	\begin{equation}\label{1.2}
		\lan y - y_0, \; x-x_0 \ran \ge 0 \quad \mbox{ for all } (x, y) \in {\rm Gr}(A).
	\end{equation}  Since $R(A+\lambda J_\vp) = X^*,$  there exists $(x_1, y_1) \in {\rm Gr}(A)$ such that
	\begin{equation}\label{1.3}
		\la J_\vp x_1 + y_1 = \la J_\vp x_0+ y_0. 
	\end{equation}
	Since $\vp$ is strictly increasing, Lemma 1.5.4 in Alber~\cite{Alber}  implies
	\begin{equation}\label{1.4}
		\lan J_\vp u_1 - J_\vp u_0, \; u_1-u_0\ran \ge \big(\vp(\|u_1\|- \vp(\|u_0\|)\big)(\|u_1\|- \|u_0\|)\ge 0
	\end{equation} for all $u_0, u_1\in X.$   
	Using $(x_1, y_1)$ in (\ref{1.2})  and simplifying with the help of (\ref{1.3}) and (\ref{1.4}), we get
	\begin{equation}\label{1.5}
		\lan J_\vp x_1 - J_\vp x_0, \; x_1-x_0\ran =0.  
	\end{equation}
	Moreover, in view of (\ref{1.4}) and (\ref{1.5}), we find
	\begin{eqnarray*}
		0&=&\lan J_\vp x_1 - J_\vp x_0, \; (x_1-x_0)/2\ran\\
		&=&\left \lan J_\vp x_1 - J_\vp \left((x_1+x_0)/2\right), \; (x_1-x_0)/2\right\ran \\
		&&+ \left \lan J_\vp \left((x_1+x_0)/2\right)-J_\vp x_0 , \; (x_1-x_0)/2\right\ran\\
		&\ge& \left[\vp(\|x_1\|)- \vp\left(\left \|(x_1+x_0)/2\right\|\right)\right]\left(\|x_1\|- \left \|(x_1+x_0)/2\right\|\right) \\
		&&+\left[\vp\left(\left \|(x_1+x_0)/2\right\|\right)- \vp(\|x_0\|)\right]\left( \left \|(x_1+x_0)/2\right\|-\|x_0\|\right) \\
		&\ge& 0.
	\end{eqnarray*}
	Consequently, $\displaystyle \|x_1\| = \|x_0\| =\left \|x_1+x_0\right\|/2 .$ Since $X$ is strictly convex, we must have $x_1 = x_0,$ and therefore $J_\vp x_1 = J_\vp x_0.$ Using this in (\ref{1.3}), we get $y_1= y_0$, i.e., $(x_0, y_0) = (x_1, y_1)\in {\rm Gr}(A),$ and hence $A$ is maximal monotone. 
	
	Conversely, suppose that $A$ is maximal monotone.  Let $\la>0$   and $\vp$ be an arbitrary gauge function.  By Lemma 1.5.5 in Alber~\cite[p. 34]{Alber},  $J_\vp$ is single-valued, strongly monotone and  demicontinuous (hence hemicontinuous) on $X$. It is clear that $J_\vp$ is coercive, i.e., $\lim_{\|x\|\to\infty} {\lan J_\vp x, x\ran }/{\|x\|} = \infty.$ Fix $y_0\in X^*$ and define 
	$Bx = \la J_\vp x- y_0$, $x\in X$.  Clearly, $B$ is also single-valued, demicontinuous, (strongly) monotone, and coercive.  By  Theorem 2.1~\cite{BA}, there exists an  $x_0\in X$ such that 
	$$\lan v+ \la J_\vp x_0- y_0, \; u -x_0\ran \ge 0 \quad \mbox{ for all } (u, v) \in {\rm Gr}(A).$$ The maximal monotonicity of $A$ implies 
	$ (x_0, -\la J_\vp x_0+y_0)\in {\rm Gr}(A), $
	i.e., 
	$y_0\in Ax_0+\la J_\vp x_0,$ and hence $A+\la J_\vp$ is surjective.
\end{proof}

The rest of the paper is organized as follows. In Section~\ref{S2}, we extend the results  in Pr\"u\ss~\cite{Pruss1981} and  Kartsatos~\cite{Kartsatos2008}  about the normalized duality mapping $J$ to   the duality mapping $J_\vp$ with an arbitrary gauge function $\vp.$  
Section~\ref{S3} concerns with  general Yosida approximants and resolvents with respect to any gauge function $\vp.$  By applying the results established in Section~\ref{S2},  we give a concise proof of the  continuity of general Yosida approximants and resolvents
on $(0, \infty) \times X.$
In Section~\ref{S4}, we indicate some areas in which the theory developed in this paper is applicable.


\section{Properties of General  Duality Mappings}\label{S2}
In this section, we  extend to an arbitrary gauge function $\vp$  the  following main result  in  Kartsatos \cite{Kartsatos2009} concerning the normalized duality mapping $J: X\to 2^{X^*},$  where $X$ is a real locally uniformly convex Banach space.\\

\noindent \textbf{Theorem A \cite{Kartsatos2009}}. \textit{Let $X$ be a  locally uniformly convex Banach space. Then  for every $R>0$ and every $x_0\in X$ there exists a nondecreasing function $\psi  = \psi (R, x_0) :\mathbf R_+\to \mathbf R_+$  such that $\psi(0) = 0$, $\psi(r)>0$ for $r>0$, and 
	$$\lan x^*- x_0^*, x-x_0\ran \ge \psi(\|x-x_0\|) \|x-x_0\|$$
	for all $x$ satisfying $\|x-x_0\|\le R$ and all $x^*\in J x$ and $x_0^*\in J x_0.$}\\

In Theorem~A, Kartsatos \cite{Kartsatos2009} actually generalized  Theorem~1 in Pr\"u\ss~\cite{Pruss1981},   which was for uniformly convex Banach spaces, to locally uniformly convex Banach spaces.  Both these authors used the normalized duality mapping $J$ with  $\vp (r) = r$.   Our main objective in this paper is to obtain Theorem~A for the duality mappings $J_\vp$  with an arbitrary gauge function $\vp$ and   use it to study the continuity of Yosida approximants and resolvents of maximal monotone operators on reflexive Banach spaces.    

In the poof of the following theorem, we follow the methodology used in \cite{Kartsatos2009} and point out the main differences that occur because of  the general duality mapping $J_\vp$. 
\begin{theorem}\label{Th3}
	Let $X$ be a  locally uniformly convex Banach space. Then  for every $R>0$ and every $x_0\in X$ there exists a nondecreasing function $\psi  = \psi (R, x_0) :\mathbf R_+\to \mathbf R_+$  such that $\psi(0) = 0$, $\psi(r)>0$ for $r>0$, and 
	$$\lan x^*- x_0^*, x-x_0\ran \ge \psi(\|x-x_0\|) \|x-x_0\|$$
	for all $x$ satisfying $\|x-x_0\|\le R$ and all $x^*\in J_\vp x$ and $x_0^*\in J_\vp x_0.$
\end{theorem}

\begin{proof}
	Fix $x_0\in X$ and define $\psi(0) =0$, $\psi(r) = \psi (R)$ for $r\ge R$,  and 
	$$\psi(r) = \inf \left\{\frac{\lan x^*-x_0^*, x-x_0\ran  }{\|x-x_0\|}\colon  x\in \overline{B_R(x_0)}\setminus B_r(x_0), \; x^*\in J_\vp x,\;  x_0^*\in J_\vp x_0 \right\},$$
	where $r\in (0, R].$ It is clear that $\psi$ is nondecreasing. We only need to prove that $\psi (r) >0$ for $r>0$.  Suppose, for a contradiction, that there exists $r_0\in (0, R]$ such that $\psi(r_0) = 0.$ Then there exists a sequence $\{x_n\}\subset \overline{ B_R(x_0)}\setminus B_r(x_0) $ and sequences  $\{x_n^*\}$ and $\{x_{0n}^*\}$ such that $x_n^*\in J_\vp x_n$,  $X_{0n}^*\in J_\vp x_0$ and
	\begin{equation}\label{2.1}
		\lim\limits_{n\to\infty}\frac{\lan x_n^*-x_{0n}^*, x_n-x_0\ran  }{\|x_n-x_0\|}=0.
	\end{equation}
	In view of Lemma 1.5.4 in \cite{Alber}, it follows from (\ref{2.1}) that
	\begin{equation*}
		\lim\limits_{n\to\infty}\Big(\vp(\|x_n\|)- \vp(\|x_0\|)\Big) (\|x_n\| - \|x_0\|) \le \lim\limits_{n\to\infty}\lan x_n^*-x_{0n}^*, x_n-x_0\ran =0.
	\end{equation*}
	for all $n.$  Since $\vp$ is strictly increasing,  we get
	$$\lim\limits_{n\to\infty}\Big(\vp(\|x_n\|)- \vp(\|x_0\|)\Big) (\|x_n\| - \|x_0\|) = 0.$$
	This implies $\|x_n\| \to \|x_0\|$  because $\{x_n\}$ is  bounded and $\vp$ is continuous and strictly increasing.  If $\|x_0\| = 0,$ then we would get $x_n\to 0;$  however, this would contradict $x_n\not\in B_r(0)$.  Therefore, $\|x_0\|>0$, and we may now assume that $\|x_n\|\ne 0$ for all $n.$  Since $X$ is locally uniformly convex, by following the arguments used in (\cite{Kartsatos2009}) we find a number  $\delta = \delta (r/\|x_0\|, x_0) >0$ such that 
	\begin{equation}\label{2.2}
		\limsup\limits_{n\to\infty}\|x_n+ x_0\| \le 2 \|x_0\| (1-\delta).
	\end{equation}
	Since  $\{\lan x_n^*, x_0\ran\}$  and $\{\lan x_{0n}^*, x_n\ran\}$  are bounded, we may assume passing to subsequences  that the limits 
	$$\lim\limits_{n\to\infty}\lan x_n^*, x_0\ran \quad\mbox{ and }\quad \lim\limits_{n\to\infty}\lan x_{0n}^*, x_n\ran $$
	exist as real numbers.
	We observe that  $$\lim\limits_{n\to\infty}\lan x_{n}^*, x_0\ran \le \limsup\limits_{n\to\infty}\vp(\|x_n\|)\|x_0\| = \vp(\|x_0\|) \|x_0\|$$
	and  $$\lim\limits_{n\to\infty}\lan x_{0n}^*, x_n\ran \le \limsup\limits_{n\to\infty}\vp(\|x_0\|)\|x_n\| = \vp(\|x_0\|) \|x_0\|.$$
	These inequalities and (\ref{2.1}) imply
	\begin{eqnarray*}
		0&\le& \Big(\vp(\|x_0\|) \|x_0\| - \lim\limits_{n\to\infty}\lan x_{n}^*, x_0\ran\Big) + \Big( \vp(\|x_0\|) \|x_0\| - \lim\limits_{n\to\infty}\lan x_{0n}^*, x_n\ran\Big)\\
		&=& \lim\limits_{n\to\infty}\Big(\vp(\|x_n\|) \|x_n\| - \lan x_n^*, x_0\ran \Big) +\lim\limits_{n\to\infty}\Big( \vp(\|x_0\|) \|x_0\| - \lan x_{0n}^*, x_n\ran \Big)\\
		&\le&\lim\limits_{n\to\infty}\lan x_n^*-x_{0n}^*, x_n-x_0\ran\\
		&=&0.
	\end{eqnarray*}
	Therefore,  $  \lim\limits_{n\to\infty}\lan x_{n}^*, x_0\ran= \lim\limits_{n\to\infty}\lan x_{0n}^*, x_n\ran = \vp(\|x_0\|) \|x_0\|,$  and
	consequently
	$$\lim \limits_{n\to\infty}\lan x_n^*, x_n +x_0\ran = \lim\limits_{n\to\infty}\lan x_n^*, x_n\ran  +\lim\limits_{n\to\infty}\lan x_n^*, x_0\ran = 2\vp(\|x_0\|) \|x_0\|.$$
	In view of (\ref{2.2}), we obtain
	\begin{eqnarray*}
		2\vp(\|x_0\|)  \|x_0\| (1-\delta)&\ge&	\vp(\|x_0\|) \limsup\limits_{n\to\infty}\|x_n+ x_0\| \\
		&= & \limsup\limits_{n\to\infty}\Big(\vp(\|x_n\|) \|x_n+ x_0\| \Big)\\
		&= & \limsup\limits_{n\to\infty}\Big(\|x_n^*\|\|x_n+ x_0\| \Big)\\
		&\ge& \lim \limits_{n\to\infty}\lan x_n^*, x_n +x_0\ran \\
		&=&2\vp(\|x_0\|) \|x_0\|,
	\end{eqnarray*}
	which  a contradiction because 	 $\|x_0\|>0.$
\end{proof}

\begin{corollary}
	Let $X$ be a  reflexive and locally uniformly convex Banach space with locally uniformly convex dual $X^*$. Then for every $R>0$ and every $x_0\in X$, there exists a function $\psi  = \psi (R, x_0) :\mathbf R_+\to \mathbf R_+$  such that $\psi(0) = 0$, $\psi(r)>0$ for $r>0$, and 
	$$\lan J_\vp x- J_\vp x_0, x-x_0\ran \ge \psi(\|x-x_0\|) \|x-x_0\|$$
	for all $x$ satisfying $\|x-x_0\|\le R.$ 
\end{corollary}

\begin{proof}
	Since $X^*$ is locally uniformly convex, it is strictly convex (see Cioranescu \cite[Proposition 2.7, p. 49]{Cioranescu}),  and therefore $X$ is smooth (see  \cite[Theorem 1.3, p. 43]{Cioranescu}). Consequently, the duality mapping $J_\vp$ is single-valued (see  \cite[Corollary 4.5, p. 27]{Cioranescu}), and  the corollary now follows immediately from Theorem~\ref{Th3}. 
\end{proof}

As a consequence of Theorem~\ref{Th3}, we directly obtain the following  corollary which is Proposition~2.17 in \cite[p. 55]{Cioranescu}.

\begin{corollary}\label{Cor} 
	Let $X$ be  a smooth and locally uniformly convex Banach space. Let $x_0\in X$ and $\{x_n\}\subset X$ be a sequence such that 
	$$\lim\limits_{n\to\infty}\lan J_\vp x_n - J_\vp x_0 , x_n - x_0\ran =0.$$ Then $x_n\to x_0$ in $X$. 
\end{corollary}

\begin{proof} In view of Lemma 1.5.4 in \cite{Alber}, we have
	\begin{equation}\label{5}
		\lim\limits_{n\to\infty}\Big(\vp(\|x_n\|)- \vp(\|x_0\|)\Big) (\|x_n\| - \|x_0\|) \le \lim\limits_{n\to\infty}\lan J_\vp x_n -J_\vp x_0, x_n-x_0\ran =0.
	\end{equation}
	This implies $\|x_n\|\to \|x_0\|.$   Therefore, there exists $R>0$ such that $x_n\in B_R(x_0)$ for all $n.$ By Theorem~\ref{Th3},  there exists a function $\psi = \psi(R, x_0)$ such that 
	$\psi(0) = 0$, $\psi(r)>0$ for $r>0$, and 
	$$\lan  J_\vp x_n- J_\vp x_0, x_n-x_0\ran \ge \psi(\|x_n-x_0\|) \|x_n-x_0\|$$
	for all $n$ satisfying $\|x-x_0\|\le R.$ This implies $x_n\to x_0$ in $X.$ 
\end{proof}

\begin{remark}\rm 
	An important class   of bounded demicontinuous operators of type $(S_+)$ is that of  the duality mappings $J_\vp$ in the settings of Corollary~\ref{Cor}.  By Lemma 1.5.5 in Alber~\cite[p. 34]{Alber},  $J_\vp$ is single-valued and  demicontinuous  on $X$. To show $J_\vp$ is of type $(S_+)$, let $\{x_n\}$ be a sequence in $X$ such that $x_n\rh x_0\in X$  and 
	$$\limsup_{n\to\infty}\lan J_\vp x_n -J_\vp x_0, x_n - x_0\ran \le 0.$$ Following the arguments in the proof of Corollary~\ref{Cor}, we get $x_n\to x_0$ in $X$. 
\end{remark}

\section{General Yosida Approximants  and Resolvents}\label{S3}
Let $X$ be a  reflexive Banach space with strictly
convex $X$ and $X^*$.
We  study variants of  the  standard Yosida approximants and resolvents  of  a maximal monotone operator $\A$, but now using the duality mapping $J_\vp$ corresponding to an arbitrary gauge function $\vp.$ 
Given  a number $\lambda>0$, a gauge function $\vp$,   and an element  $x\in X$,  the inclusion 
\begin{equation} \label{3.1}
	0\in J_\vp(x_\lambda- x) + \lambda Ax_\lambda 
\end{equation} has a unique solution $x_\lambda\in D(A)$.     We point out that the existence of a solution $x_\lambda$  of (\ref{3.1}) is ensured by using  $J_\vp (\cdot - x)$  in place of $J_\vp$ in Theorem~\ref{Th1},  and the uniqueness of the solution $x_\la$ follows from the strict convexity of $X$.
We now define $J_\lambda^\vp : X\to D(A)\subset  X$  and $A_\lambda^\vp : X\to X^*$ by 
\begin{equation}\label{3.2}
	J_\lambda^\vp x := x_\lambda \quad \mbox{ and }\quad A_\lambda^\vp x: = \dfrac{1}{\lambda} J_\vp (x- J_\lambda^\vp x), \quad x\in X.
\end{equation} The operators $A_\lambda^\vp$ and $J_\lambda^\vp$ corresponding to an arbitrary gauge function $\vp$  are generalizations of the  standard Yosida approximant $A_\lambda$  and resolvent $J_\lambda$  of $A$ , respectively. 
We see that \begin{equation}\label{3.3}
	A_\lambda^\vp x\in A (J_\lambda^\vp x) \quad \mbox{and} \quad 	x = J_\lambda^\vp x + J^{-1}_\vp ( \la A_\lambda^\vp x).  
\end{equation}
It is easy to verify  that $A= A_\lambda^\vp $ if and only if $A=0$. In fact, if $A= 0$, then $J_\lambda^\vp   = I$, the identity operator on $X$. Moreover,  $0\in D(A)$ and $0\in A(0)$ if and only if  $A_\lambda ^\vp 0 =0$.	 

When $\vp(r) = r,$ a splitting of $x$  in terms of $A_\lambda$ and $J_\lambda$ is 
\begin{equation*}
	x = J_\lambda x + \lambda J^{-1} ( A_\lambda x),
\end{equation*}  and consequently
\begin{equation}\label{3.4}
	A_\lambda  x= \left(A^{-1}+ \lambda  J^{-1}\right)^{-1} x, \quad x\in X.
\end{equation}
Such  a closed form for $A_\lambda^\vp$ is also available for $\vp(r) = r^{p-1}, \; p\in (1, \infty);$  but not for a general $\vp$.

A proof of the following lemma  for the case $\vp(r) = r$  can be found in Boubakari and Kartsatos \cite{BK}.  Recently in \cite{Adhikari2022},  Adhikari et al.  proved the following lemma for an arbitrary gauge function $\vp$. 
\begin{lemma}\label{L3}
	Let $\A$ be  maximal monotone  and $G\subset X$ be  bounded. Let  $0<\lambda_1<\lambda_2$. Then there exists  a constant $K$, independent of $\lambda$,  such that 
	$$\|A_\lambda^\vp x\| \le K$$
	for all $x\in \overline{ G}$ and $\lambda\in [\lambda_1, \lambda_2].$
\end{lemma}

By a well-known renorming theorem due to Trojanski \cite{Trojanski},
one can always renorm a reflexive Banach space $X$ with an equivalent
norm with respect to which both $X$ and $X^*$ are still dual to each other and both become locally uniformly convex (therefore strictly convex).  After such a renorming, the  duality mapping $J_\varphi $ is single-valued homeomorphism from $X$ onto $X^*$, and the inverse $J_\vp^{-1}$ is the duality mapping $ J_{\vp^{-1}}$ of $X^*$ with the gauge function $\vp^{-1}$.    Henceforth we assume that $X$ is reflexive and both $X$ and $X^*$ are locally uniformly convex.	

For the special case  $\vp(r) = r$, the reader can find a proof of Theorem~\ref{Th2} below  in Kartsatos and Skrypnik \cite{KartsatosSkrypnik2005a}  when $0\in A(0)$ and in  Asfaw and Kartsatos \cite{ASK2012},   without the condition $0\in A(0)$.  We note that Zhang and Chen in \cite[Lemma~2.7, p. 446]{ZC} proved the continuity of $x\mapsto A_\lambda x$   on $ D(A)$ for each $\lambda>0$, also  without the condition $0\in A(0)$.   In  the proof of \cite[Lemma~6]{ASK2012}, however,   the continuity of $x\mapsto A_\lambda x$ on $X$ is  used without a justification .  In a recent paper of Adhikari et al., there is a proof for $\vp(r) = r^{p-1}$, $p\in (0, \infty),$  in which the authors used the explicit form (\ref{3.4}) for $A_\la^\vp$  just as in proofs of  the continuity of $ (\la, x)\mapsto  A_\la x$  on $(0, \infty)\times X$ in both \cite{ASK2012} and  \cite{KartsatosSkrypnik2005a}.   However,  the explicit form (\ref{3.4})  is not available for $A_\la^\vp$ when $\vp$ is  nonhomogeneous. Our  method  needs no such explicit form, and we next provide  a concise proof of the  continuity of the mappings $(\lambda, x)\mapsto A_\lambda^\vp x$  and $(\lambda, x)\mapsto J_\lambda^\vp x$ on $(0, \infty) \times X$ for an arbitrary gauge function $\vp$.

\begin{theorem}\label{Th2}
	Let $\A$ be  a maximal monotone operator and $\vp$ be an arbitrary gauge function.  Then the mappings $(\lambda, x)\mapsto A_\lambda^\vp x$ and $(\lambda, x)\mapsto J_\lambda^\vp x$ are continuous on $(0, \infty) \times X$. 
\end{theorem}
\proof  It suffices to prove the continuity of  $(\lambda, x)\mapsto A_\lambda^\vp x$ on $(0, \infty) \times X$ because the continuity of $(\lambda, x)\mapsto J_\lambda^\vp x$ on $(0, \infty) \times X$  follows from (\ref{3.2}) and the continuity of $J_\vp^{-1}.$  To this end, let $\{\lambda_n\}\subset (0, \infty)$ and $\{x_n\}\subset X$ be such that $\lambda_n\to \lambda_0\in (0, \infty)$ and $x_n\to x_0\in X$. Let $G\subset X$ be a bounded set that contains $x_n$ for all $n$.  Rename $\lambda_1, \lambda_2>0$ such that $\lambda_n\in [\lambda_1, \lambda_2]$ for all $n$.  Denote $u_n=J_{\lambda_n}^\vp x_n$ and  $v_n = A_{\lambda_n}^\vp x_n$ for all $n = 0, 1, 2, ....$  We prove that $v_n \to v_0$. 
It follows from Lemma~\ref{L3} that the sequence $\{v_n\}$ is bounded,  and therefore
\begin{equation}\label{3.5}
	\lim\limits_{n\to\infty}\lan \la_n v_n- \la_0v_0, x_n-x_0\ran = 0.
\end{equation}
In view of (\ref{3.3}), we get
\begin{equation}\label{3.6}
	x_n =u_n + J_\vp^{-1}(\la_nv_n)  \mbox{ and } x_0 =u_0 + J_\vp^{-1}(\la_0v_0).
\end{equation}
Using (\ref{3.6}), we have
\begin{equation}\label{3.7}
	\begin{split}
		\lan\la_n &v_n- \la_0v_0,\; x_n-x_0\ran\\=&\;\lan \la_nv_n- \la_0v_0,\;  u_n-u_0\ran\\&+ \lan \la_nv_n- \la_0v_0, \; J_\vp^{-1}(\la_nv_n)- J_\vp^{-1}(\la_0v_0)\ran\\
		=&\;\la_n\lan v_n- v_0,\;  u_n-u_0\ran +\lan (\la_n- \la_0)v_0,\;  u_n-u_0\ran\\
		&+ \lan \la_nv_n- \la_0v_0, \; J_\vp^{-1}(\la_nv_n)- J_\vp^{-1}(\la_0v_0)\ran.
	\end{split}
\end{equation}
Since $\{\la_n\}$, $\{x_n\}$ and $\{v_n\}$ are bounded and $J_\vp^{-1}$ is bounded,   it follows from (\ref{3.6}) that $\{u_n\}$ is bounded. 
Since  $A$ is monotone and $\la_n\to \la_0$, it follows from (\ref{3.5}) and (\ref{3.7}) that 
\begin{equation*}
	\lim\limits_{n\to\infty}\lan \la_nv_n- \la_0v_0, J_\vp^{-1}(\la_nv_n)- J_\vp^{-1}(\la_0v_0)\ran = 0.
\end{equation*}
Since $J_\vp^{-1}$ is a duality mapping of $X^*,$
applying Corollary~\ref{Cor} to $J_\vp^{-1}$ gives
$\la_nv_n\to \la_0 v_0$ in $X^*.$  Since $\la_n\to\la_0>0$, we obtain $v_n\to  v_0$ in $X^*,$  which proves the continuity  of $A_\la^\vp$ on $(0, \infty)\times X$.
\endproof


\section{Discussion}\label{S4}  An important reason for the establishment of the continuity of $A_\la^\vp$ is to make sure that  Browder's degree in \cite{Browder1983} for $A+F$, where $\A$ is  maximal monotone,  and $F: \overline G\to X^*$ is bounded demicontinuous operators of type $(S_+)$ with $G\subset X$ is open and bounded,  is unaffected even when the normalized duality mapping $J$  used in the  Browder theory is replaced with the duality mapping $J_\vp$ for any gauge function $\vp$.      Let $\la_1, \la_2>0$ and define $q(t) := \la_1t+ (1-t) \la_2$,  $t\in [0, 1]$. Using the continuity of $A_\la^\vp$ established in Theorem~\ref{Th3},  we  verify that $T_t:=A_{q(t)}^\vp $ is a \textit{pseudomonotone homotopy}  in the sense of Browder~\cite[Definition 8, p. 32]{Browder1983}; \textit{i.e, 
 for a given sequence $\{t_n\}\subset [0, 1]$ converging to $t_0\in [0, 1]$ and 	 a given $(x, y) \in {\rm Gr}(T_{t_0}),$ there exists a 
	sequence $\{(x_n, y_n)\}$ such that $(x_n, y_n) \in{\rm Gr}(T_{t_n}),$   $x_n\to  x_0$ in $X$ and  $y_n\to  y_0$ in $X^*.$ }
In fact,   let $\{t_n\}\subset [0, 1]$ be such that $t_n\to t_0\in [0, 1]$ and let $(x_0, y_0) \in {\rm Gr} (T_{t_0}).$   Define $x_n = x_0$ and $y_n = T_{t_n} x_n = A_{q(t_n)}^\vp x$ for all $n$.  It is clear that $x_n \to x_0$  in $X$ and, by Theorem~\ref{Th3}, $y_n \to A_{q(t_0)}^\vp x_0 = T_{t_0}x_0 =y_0$ in $X^*$. 

Berkovits and Miettunen in \cite{BM2008} established the uniqueness of  Browder's degree in \cite{Browder1983}, and therefore  Browder's original degree defined  in terms of the degree, ${\rm d}_{(S_+)}$, of bounded, demicontinuous  and $(S_+)$-mappings by
$$\lim\limits_{\la\to 0^+} {\rm d}_{(S_+)}(A_\la+F, G, p) $$ coincides with $$\lim\limits_{\la\to 0^+} {\rm d}_{(S_+)}(A_{\la }^\vp +F, G, p),$$ 
where $p\in X^*$ satisfies $p\not\in (A+F)(\partial G\cap D(A)).$ We  point out that Kien et al. \cite{Kien2008}  used the definition of pseudomonotone homotopy as given above and enlarged the class of  pseudomonotone homotopies by not requiring $0\in T_t(0)$ for all $t$.

\begin{definition}\label{Hom} \rm
	An operator $A: X\supset D(A)\to 2^{X^*}$ is said to be
	\textit{positively homogeneous} of degree $\gamma >0$ if $(x,y)\in {\rm Gr}(A)$ implies $sx\in D(A)$ for all
	$s\ge 0$ and $(sx, s^\gamma y)\in {\rm Gr}(A).$ 
\end{definition}

Since $J$ is homogeneous of degree $1$,  it follows from (\ref{1.1}) that
\begin{equation*}
	J_\varphi(s x) = \frac{\vp(s\|x\|)}{\|x\|} Jx, \quad (s, x)\in \mathbf{R}_+\times X.
\end{equation*}
In particular, when $\vp(r) = r^{p-1}, 1< p< \infty$, we get  	$J_\varphi x  = \|x\|^{p-2} Jx,$  $x\in X,$ which implies
\begin{equation*}
	J_\varphi(s x) = s^{p-1} J_\vp x, \quad (s, x)\in \mathbf{R}_+\times X.
\end{equation*}
Thus, $J_\varphi$ is positively homogeneous of degree $p-1.$ 
When $X$ and $X^*$ are reflexive and strictly convex,  the inverse $J^{-1}_\vp$ of the single-valued $J_\vp$ is the duality mapping $J_{\vp^{-1}}$ of $X^*$ with  $\vp^{-1}(r) = r^{q-1},$ where $q$ satisfies $1/p + 1/q = 1.$  Consequently,
\begin{equation*}
	J_{\varphi^{-1}}(s x^*) = s^{q-1} J_{\varphi^{-1} }(x^*), \quad (s, x^*)\in \mathbf{R}_+\times X^*.
\end{equation*}
In general, $J_\vp$ is positively homogeneous of degree $\gamma>0$ if and only if $\vp$ is positively homogeneous of  degree $\gamma>0.$ 

The duality mappings $J_\vp$  for the gauge functions $\vp(r) = r^{p-1}, \; p\in (1, \infty),$ play an important role in  nonlinear Fredholm theory concerning maximal monotone operators that are   positively homogeneous of degree exactly $p-1$.  In fact, for  $\vp(r) = r^{p-1},$ $p\in (1, \infty),$   the Yosida approximants $A_\la^\vp$ of a positively homogeneous maximal monotone operator $A$ of degree $p-1$  are also positively homogeneous of the same degree, among other properties. This result is established in a recent work of Adhikari et al. \cite{Adhikari2022}.  For more information about the existence theorems of  Fredholm alternative type that involve positively homogeneous maximal monotone operators, the reader is referred to Hess~\cite{Hess}.

\end{document}